\documentclass{article}

\usepackage{natbib}            
\usepackage{graphicx}          

\usepackage{amsmath}
\usepackage{amssymb}
\usepackage{amsthm}

\usepackage{maplestd2e}

\newcommand{\Z}{\mathbb{Z}}

\newcommand{\C}{\mathbb{C}}

\newcommand{\ld}{{\rm ld}}
\newcommand{\init}{{\rm init}}
\newcommand{\disc}{{\rm disc}}
\newcommand{\sep}{{\rm sep}}

\newcommand{\lcm}{{\rm lcm}}

\newcommand{\Sol}{{\rm Sol}}

\newcommand{\NF}{{\rm NF}}

\newtheorem{thm}{Theorem}
\newtheorem{prop}[thm]{Proposition}
\newtheorem{cor}[thm]{Corollary}

\theoremstyle{definition}
\newtheorem{defn}[thm]{Definition}
\newtheorem{rem}[thm]{Remark}
\newtheorem{exmp}[thm]{Example}
\newtheorem*{ack}{Acknowledgement}

\begin{document}

\title{Thomas decompositions of parametric nonlinear control systems} 


\author{Markus Lange-Hegermann} 
\author{Daniel Robertz}

\begin{abstract}                          
This paper presents an algorithmic method to study
structural properties of nonlinear control systems in
dependence of parameters. The result consists of a
description of parameter configurations which cause
different control-theoretic behaviour of the system
(in terms of observability, flatness, etc.).
The constructive symbolic method is based on the differential
Thomas decomposition into disjoint simple systems,
in particular its elimination properties.
\end{abstract}

\section{Introduction}\label{sec:intro}

Symbolic computation allows to study many structural aspects
of control systems, e.g., controllability, observability,
input-output behaviour, etc. In contrast to a numerical
treatment, the dependence of the results on parameters occurring
in the system is accessible to symbolic methods.

An algebraic approach for treating nonlinear control systems
has been developed during the last decades,
cf., e.g., \cite{FliessGlad},
\cite{DiopElimination,DiopDiffAlg}, \cite{Pommaret2001},
and the references therein. In particular, the notion of
flatness has been studied extensively and has been applied to many
interesting control problems
(cf., e.g., \cite{FliessLevineMartinRouchon}).
The approach of Diop builds on the characteristic set method
(cf.\ \cite{Kolchin}, \cite{Wu}). The Rosenfeld-Gr{\"o}bner algorithm
(cf.\
\cite{BoulierLazardOllivierPetitotAAECC})
can be used to perform the relevant computations effectively;
cf.\ also \cite{WangEliminationMethods} for alternative approaches.

Dependence of control systems on parameters has been examined, in
particular, in \cite{PommaretQuadrat1997}, \cite{Pommaret2001}.
In the case of linear systems, stratifications of the space of
parameter values have been studied using Gr{\"o}bner
bases in \cite{LevandovskyyZerz}.

In the 1930s the American mathematician J.~M.\ Thomas designed an
algorithm which decomposes a polynomially nonlinear system of
partial differential equations into so-called simple systems.
The algorithm uses, in contrast to the characteristic set method,
inequations to provide a disjoint decomposition of the solution
set (cf.\ \cite{Thomas}).
It precedes \cite{Kolchin} and \cite{Seidenberg} (building on
\cite{Ritt}).
Recently a new algorithmic approach to
the Thomas decomposition method has been developed
(cf.\
\cite{GerdtThomas,BaechlerGerdtLangeHegermannRobertz2,Robertz2012}),
also using ideas of M.~Janet (cf.\ \cite{Janet1});
cf.\ also \cite{WangSimpleSystems} for an earlier implementation
of the algebraic part.

So far the dependence of nonlinear control systems on parameters has not
been studied by such a rigorous method as the Thomas decomposition.
This paper demonstrates how the Thomas decomposition method can be
applied in this context. In particular, the Thomas decomposition can detect certain
structural properties of control systems by performing elimination
and it can separate singular cases of behaviour in control systems from
the generic case due to splitting into disjoint solution sets.

In Section~\ref{sec:thomas} we sketch the Thomas decomposition
method for algebraic and differential systems. The algorithm for
the differential case builds on the algebraic part.
Section~\ref{sec:elimination} describes how the differential
Thomas decomposition can be used to solve elimination problems
that occur in our study of nonlinear control systems.
In Section~\ref{sec:systemstheory} we recall the notions of
observability and flatness, which are addressed in the examples
in Section~\ref{sec:applications} using our Maple
implementation.

\section{Thomas Decomposition}\label{sec:thomas}

\subsection{Simple Algebraic Systems}\label{subsec:algsimple}

Let $K$ be a field of characteristic zero. We denote by
$R$ the polynomial algebra $K[x_1, \ldots, x_n]$ and
we fix a total ordering $<$ on $\{ x_1, \ldots, x_n \}$.
Then the greatest variable with respect to $<$ which
occurs in a non-constant polynomial $p \in R$ is called
the \emph{leader} of $p$ and is denoted by $\ld(p)$.

In what follows, we regard every $p \in R \setminus K$ as
a polynomial in $\ld(p)$ with coefficients in
$K[x_i \mid 1 \le i \le n, \, x_i < \ld(p)]$, and, recursively,
each coefficient as a polynomial in its leader, etc.
The coefficient of the highest power of $\ld(p)$ in $p$
is called the \emph{initial} of $p$, denoted by $\init(p)$.
Moreover, we denote the
discriminant of $p$ (as a polynomial in $\ld(p)$) by
$\disc(p)$. Both
$\init(p)$ and $\disc(p)$ are elements of the polynomial ring
$K[x_i \mid 1 \le i \le n, \, x_i < \ld(p)]$.

Let $\overline{K}$ be an algebraic closure of $K$. In this
subsection we consider
(finite) algebraic systems (of equations and inequations, defined over $K$)
\begin{equation}\label{eq:algebraicsystem}\tag{A}
S = \{ p_1 = 0, \, \ldots, \, p_s = 0, \, q_1 \neq 0, \, \ldots, \, q_t \neq 0 \},
\end{equation}
where $p_i$, $q_j \in R$, $s$, $t \in \Z_{\ge 0}$.
The set of solutions of $S$ (in $\overline{K}^n$) is defined as
\[
\Sol_{\overline{K}}(S) = \{ a \in \overline{K}^n \mid p_i(a) = 0, \, q_j(a) \neq 0 \mbox{ for all } i, j \}.
\]
Finally, for $0 \le k < n$,
let $\pi_k : \overline{K}^n \to \overline{K}^{n-k}$ be the projection
onto the coordinate subspace of dimension $n-k$ whose coordinates are ranked
lowest with respect to $<$, i.e.,
$\pi_k(a_1, \ldots, a_n) = (a_{i_j} \mid 1 \le j \le n-k)$,
where $1 \le i_1 < \ldots < i_{n-k} \le n$ and
$x_{i_1}$, \ldots, $x_{i_{n-k}}$ are the smallest $n-k$ elements
of $\{ x_1, \ldots, x_n \}$ with respect to $<$.

\smallskip

\begin{defn}\label{de:algsimple}
An algebraic system $S$ as in (\ref{eq:algebraicsystem})
is \emph{simple} (with respect to $<$)
if the following conditions hold.
\begin{enumerate}
\item $p_i \not\in K$, $q_j \not\in K$ for all $i$ and $j$.
\item $|\{ \ld(p_1), \ldots, \ld(p_s), \ld(q_1), \ldots, \ld(q_t) \}| = s+t$.
\item For every $r \in \{ p_1, \ldots, p_s, q_1, \ldots, q_t \}$,
if $\ld(r)$ is the $k$-th greatest variable w.r.t.\ $<$, then neither $\init(r) = 0$ nor
$\disc(r) = 0$ has a solution in $\pi_k(\Sol_{\overline{K}}(S))$.
\end{enumerate}
\end{defn}

\smallskip

\begin{rem}
Condition (2) implies that the leaders of the left hand sides
of the equations and inequations in a simple algebraic system
$S$ are pairwise distinct, i.e., $S$ is a triangular set
(cf., e.g., \cite{AubryLazardMorenoMaza,HubertIandII,WangEliminationMethods};
cf.\ \cite{RegularChains} for a related implementation).

The meaning of condition (3) is that for every ${1 \le k < n}$ and
each $(a_1, \ldots, a_{n-k}) \in \pi_k(\Sol_{\overline{K}}(S))$
there exists ${a \in \overline{K}}$ such that
${(a_1, \ldots, a_{n-k}, a) \in \pi_{k-1}(\Sol_{\overline{K}}(S))}$,
and the number of possible values for $a$ is either finite or co-finite,
determined by the degree of the leader in the corresponding
equation or inequation in $S$, if any
(cf.\ also \cite{PleskenCounting}).
\end{rem}

\smallskip

\begin{defn}\label{de:algthomas}
Let $S$ be an algebraic system as in (\ref{eq:algebraicsystem}).
A \emph{Thomas decomposition} of $S$ (with respect to $<$)
consists of finitely many
simple algebraic systems ${S_1, \ldots, S_k}$ such that
$\Sol_{\overline{K}}(S)$ is the disjoint union of
$\Sol_{\overline{K}}(S_1)$, \ldots, $\Sol_{\overline{K}}(S_k)$.
\end{defn}

\smallskip

\begin{rem}\label{rem:algebraicthomasalg}
Euclidean pseudo-division and splitting of systems
allow to compute a Thomas decomposition for any given
(finite) algebraic system (defined over a computable field $K$
of characteristic zero) in finitely many steps.

In each round the algorithm chooses a system $S$ from a set of
systems to be processed. If $S$ contains an equation
whose left hand side is a non-zero constant or an inequation
with zero left hand side, then $S$ is discarded because its
solution set is empty. Moreover, we assume that equations
$0 = 0$ and inequations $k \neq 0$, where $k \in K \setminus \{ 0 \}$
are removed from each system, so that now the leader of every
equation and inequation is well-defined.

The algorithm applies Euclidean pseudo-division to each
pair $p_1 = 0$, $p_2 = 0$ of distinct equations in $S$ with
$\ld(p_1) = \ld(p_2) =: x$; i.e., if $\deg_{x}(p_1) \ge \deg_{x}(p_2)$,
then polynomial division is performed on $c \cdot p_1$ and
$p_2$, where $c$ is a suitable power of $\init(p_2)$.
In order not to change the solution set, when $p_1 = 0$ is
replaced with the result of the polynomial division, non-vanishing
of $\init(p_2)$ on the set of solutions of $S$ is assumed. To
this end, a preparatory step splits a system into two and adds
$\init(p_2) \neq 0$ and $\init(p_2) = 0$, respectively, if
necessary.

Each pair $q_1 \neq 0$, $q_2 \neq 0$ of distinct inequations
in $S$ with $\ld(q_1) = \ld(q_2)$ is replaced with
$\lcm(q_1, q_2) \neq 0$.
The computation of the least common multiple depends on case
splittings according to vanishing of initials, similar to
the previous remarks.

For each pair $p = 0$, $q \neq 0$ in $S$ with
$\ld(p) = \ld(q)$, the algorithm computes
$r := \gcd(p, q)$ using Euclidean pseudo-division.
This process again relies on non-vanishing initials, which is
ensured by distinguishing cases appropriately in advance.
If $p$ divides $q$, then $S$ has no solution and is discarded.
If $r$ is a non-zero constant, then $q \neq 0$ is removed from $S$.
Else, $p = 0$ and $q \neq 0$ are replaced with $p/r = 0$ and
$q / r \neq 0$, respectively.

Finally, there is some flexibility when to take care of
non-vanishing discriminants. Since $K$ is of characteristic zero,
this property is established for any equation $r = 0$ or
inequation $r \neq 0$ by examining
$g := \gcd(r, \frac{\partial r}{\partial \ld(r)})$, again with
appropriate case distinctions, and thus determining the
square-free part $r / g$.
\end{rem}

We refer to \cite{BaechlerGerdtLangeHegermannRobertz2,PleskenCounting,Robertz2012}
for more information and to \cite{ThomasPackages} for an
implementation in Maple. In practice, we apply subresultants
for the computation of gcds and the related case distinctions.

\subsection{Simple Differential Systems}\label{subsec:diffsimple}

Let $K$ be a differential field of characteristic zero with
pairwise commuting derivations $\delta_1$, \ldots, $\delta_n$,
i.e., $K$ is a field of characteristic zero, and each
$\delta_i$ is a map $K \to K$ satisfying
$\delta_i(k_1 + k_2) = \delta_i(k_1) + \delta_i(k_2)$
and the Leibniz rule
$\delta_i(k_1 \, k_2) = k_1 \, \delta_i(k_2) + \delta_i(k_1) \, k_2$
for all $k_1$, $k_2 \in K$.

In this subsection we denote by $R$ the \emph{differential polynomial ring}
$K\{ u_1, \ldots, u_m \}$ in the differential indeterminates
$u_1$, \ldots, $u_m$ with
pairwise commuting derivations $\partial_1$, \ldots, $\partial_n$
whose restrictions to $K$ are $\delta_1$, \ldots, $\delta_n$, i.e.,
$R = K[(u_i)_J \mid 1 \le i \le m, \, J \in (\Z_{\ge 0})^n]$ is the
polynomial ring in the algebraically independent variables
$(u_i)_J$, and $\partial_j : R \to R$ is the derivation defined by
extending $\partial_j((u_i)_J) := (u_i)_{J+1_j}$ additively to $R$
such that it satisfies the Leibniz rule on $R$ and restricts to
$\delta_j$ on $K$. Here $1_j$ is the multi-index of length $n$
whose $i$-th entry is $1$ if $i = j$ and $0$ otherwise.

Each variable $(u_i)_J$
is thought of as representing the partial derivative, corresponding to
the multi-index $J$,
of a smooth (or rather analytic) $K$-valued function of $n$ arguments.
When dealing with differential polynomials $p \in R$
algorithmically, we compare terms with respect to a ranking, which
is defined next. We agree that applying derivations should make terms
larger with respect to the ranking. This is taken
into account as follows.

\smallskip

\begin{defn}\label{de:ranking}
A \emph{ranking} on 
$R$ is a total ordering $<$ on
$\{ (u_i)_J \mid 1 \le i \le m, \, J \in (\Z_{\ge 0})^n \}$
such that $u_i < \partial_j u_i$ for all $i$ and $j$, and such that
$(u_{i_1})_{J_1} < (u_{i_2})_{J_2}$ implies
$\partial_j (u_{i_1})_{J_1} < \partial_j (u_{i_2})_{J_2}$ for all
$j$, $i_1$, $i_2$, $J_1$, $J_2$.
\end{defn}

Every ranking is a well-ordering, i.e., there exist no infinitely
decreasing sequences of variables $(u_i)_J$.

We define
$\partial^J := \partial^{J_1} \ldots \partial^{J_n}$, $J \in (\Z_{\ge 0})^n$,
and write ${|J| := J_1 + \ldots + J_n}$ for the length of the
multi-index $J$.

\smallskip

\begin{exmp}\label{ex:degrevlex}
The \emph{degree-reverse lexicographical ranking} on the
differential polynomial ring $K\{ u \}$ (i.e., $m = 1$)
with pairwise commuting derivations
$\partial_1$, \ldots, $\partial_n$ is defined by
\[
u_{J} < u_{J'} :\Leftrightarrow \left\{
\begin{array}{l}
|J| < |J'| \mbox{ or }
(|J| \! = \! |J'| \mbox{ and } J \neq J' \mbox{ and } \\[0.2em]
J_i > J'_i \mbox{ for } i \! = \! \max \{ 1 \le k \le n \mid J_k \neq J'_k
\}).
\end{array} \right.
\]
In this example, we have
$\partial_n u < \partial_{n-1} u < \ldots < \partial_1 u$.
There are in fact $n!$ different degree-reverse lexicographical rankings
according to the ordering of the $\partial_i u$.
\end{exmp}

In this subsection we consider (finite) differential systems
(i.e.\ systems of differential equations and inequations)
\begin{equation}\label{eq:differentialsystem}\tag{D}
S = \{ p_1 = 0, \, \ldots, \, p_s = 0, \, q_1 \neq 0, \, \ldots, \, q_t \neq 0 \},
\end{equation}
where $p_i$, $q_j \in R$, $s$, $t \in \Z_{\ge 0}$.

In what follows, we assume that a ranking $<$ on $R$ is fixed.
Carrying over the concepts of Subsection~\ref{subsec:algsimple},
the \emph{leader} $\ld(p)$ of
$p \in R \setminus K$ is defined to be the
greatest variable with respect to $<$ which occurs in $p$. The
\emph{initial} of $p$, denoted by $\init(p)$, is the coefficient of the
highest power of $\ld(p)$ in $p$, and the \emph{separant} of $p$,
denoted by $\sep(p)$, is the formal partial derivative of $p$
with respect to $\ld(p)$. The initial of $p$ is an
element of the polynomial ring
$K[(u_i)_J \mid 1 \le i \le m, \, J \in (\Z_{\ge 0})^n, \, (u_i)_J < \ld(p)]$.

In order to ensure formal integrability for the kind of
differential systems we are heading for, we apply the
concept of Janet division (cf.\ \cite{Janet1,GerdtBlinkov98a}), which
restricts the usual divisibility relation on the free
commutative monoid $\Theta$ generated by $\partial_1$, \ldots, $\partial_n$.

\smallskip

\begin{defn}\label{de:admissiblederiv}
Let $M \subset (\Z_{\ge 0})^n$ be finite. \emph{Janet division}
associates with each $I \in M$ a partition of
$\{ \partial_1, \ldots, \partial_n \}$ into the subsets of
\emph{admissible} and \emph{non-admissible derivations} as follows.
Let $I = (I_1, \ldots, I_n) \in M$. Then $\partial_i$ is
admissible for $I$ if and only if
\[
I_i = \max \{ I'_i \mid (I'_1, \ldots, I'_n) \in M, \, I'_j = I_j
\mbox{ for all } 1 \le j < i \};
\]
otherwise $\partial_i$ is non-admissible for $I$. We denote by
$\mu(I, M)$ the set of derivations that are admissible for $I$ (with
respect to $M$),
$\overline{\mu}(I, M) := \{ \partial_1, \ldots, \partial_n \} \setminus \mu(I, M)$,
and
\[
\Theta(I, M) :=
\{ \partial^J \mid
J \in (\Z_{\ge 0})^n, \,
J_j = 0 \mbox{ if } \partial_j \in \overline{\mu}(I, M) \}.
\]
\end{defn}

Janet division defines for each $\partial^{I}$ with $I \in M$
a cone of multiples $\Theta(I, M) \partial^{I}$ such that each
two cones are disjoint. By enlarging $M$ appropriately, these
cones cover the set of all multiples of $\partial^{I}$, $I \in M$.

\smallskip

\begin{defn}\label{de:janetcomplete}
A finite subset $M$ of $(\Z_{\ge 0})^n$ is said to be
\emph{Janet-complete} if
$\bigcup_{I \in M} \Theta \partial^{I} = \bigcup_{I \in M} \Theta(I, M) \partial^{I}$.
\end{defn}

\smallskip

\begin{defn}\label{de:reduced}
Let $p_1$, \ldots, $p_r \in R \setminus K$,
$\ld(p_j) = \partial^{I_j} u_{i_j}$, be such that
each $M_k = \{ I_j \mid 1 \le j \le r, \, i_j = k \}$
is Janet-complete.
A differential polynomial $p \in R$ is \emph{Janet-reduced}
modulo $p_1$, \ldots, $p_r$ if, for all $1 \le j \le r$,
$\deg_{x}(p) < \deg_{x}(p_j)$, where $x := \ld(p_j)$, and no variable
$\partial^J \ld(p_j)$, where $\partial^J \in \Theta(I_j, M_{i_j})$,
$|J| > 0$, occurs in $p$.
\end{defn}

\smallskip

\begin{rem}\label{rem:janetreduction}
Any $p \in R$ can be transformed into
a differential polynomial which is Janet-reduced modulo
$p_1$, \ldots, $p_r$ by applying
Euclidean pseudo-division modulo $p_1, \ldots, p_r$ and
their derivatives repeatedly.
The obvious strategy is to first eliminate
the greatest variable $\partial^J \ld(p_j)$ with respect to $<$
that occurs in $p$ and to proceed to lower variables. If
$\partial^J \ld(p_j)$ is a proper derivative of $\ld(p_j)$, then,
before reducing modulo $p_j$,
the pseudo-division multiplies $p$ with $\sep(p_j)$
(which, in fact, is the coefficient of the leader of any proper
derivative of $p_j$). If $|J| = 0$, then the pseudo-division
multiplies $p$ with $\init(p_j)$.
As in Subsection~\ref{subsec:algsimple}, cases of vanishing
initials or separants have to be examined separately, which
results into splittings of differential systems. Under the
assumption of non-vanishing initials and separants we denote
the result of the pseudo-reduction by $\NF(p, \{ p_1, \ldots, p_r \})$.
\end{rem}

\smallskip

\begin{defn}\label{de:passive}
In the situation of Definition~\ref{de:reduced},
$\{ p_1, \ldots, p_r \}$ is \emph{passive} if, for all
$1 \le j \le r$,
$\NF(d \, p_j, \{ p_1, \ldots, p_r \}) = 0$ for all
$d \in \overline{\mu}(I_j, M_{i_j})$.
\end{defn}

\smallskip

\begin{defn}
We call a differential system $S$ as in
(\ref{eq:differentialsystem}) \emph{simple} (with respect to $<$)
if the following conditions hold.
\begin{enumerate}
\item $S$ is simple as an algebraic system (in the finitely many
variables $(u_i)_J$ which occur in $S$, totally ordered by the ranking $<$).
\item $\{ p_1, \ldots, p_s \}$ is passive.
\item $q_1$, \ldots, $q_t$ are Janet-reduced modulo $p_1$, \ldots, $p_s$.
\end{enumerate}
\end{defn}

From now on we assume that $K$ is a differential field of complex
meromorphic functions on a connected open subset $\Omega$ of $\C^n$
with coordinates $z_1$, \ldots, $z_n$, and that the
derivations $\delta_1$, \ldots, $\delta_n$ are defined
by partial differentiation with respect to $z_1$, \ldots, $z_n$,
respectively. Differential equations $p = 0$ (and their derivatives)
translate into algebraic equations
for the Taylor coefficients of a power series expansion of
$u_1$, \ldots, $u_m$ around an arbitrary point in $\Omega$.
A differential inequation $q \neq 0$ can be interpreted as the
disjunction of algebraic inequations for all Taylor coefficients
of the analytic expression that is obtained from $q$ by
substitution of power series expansions for $u_1$, \ldots, $u_m$.

In considering complex analytic functions on $\Omega$ as solutions
to differential systems, we 
assume that $\Omega$ is chosen appropriately with regard to
given differential systems, more precisely, that for every
problem instance treated by the Thomas decomposition method,
for each of the resulting simple differential systems
there exists an analytic solution on $\Omega$.
Which subsets $\Omega$ of $\C^n$ are appropriate
can often be decided only after the formal treatment of the
given differential systems by the Thomas algorithm. Questions
concerning the radius of convergence of analytic solutions
are ignored here.

For any differential system $S$ as in (\ref{eq:differentialsystem})
we denote by $\Sol_{\Omega}(S)$ the
set of complex analytic functions $f: \Omega \to \C$
satisfying $p_i(f) = 0$ and $q_j(f) \neq 0$ for all $i$ and $j$.

\smallskip

\begin{defn}
Let $S$ be a differential system as in (\ref{eq:differentialsystem}).
We call a family of finitely many simple differential systems 
$S_1, \ldots, S_k$ such that $\Sol_{\Omega}(S)$ is the disjoint 
union of the solution sets $\Sol_{\Omega}(S_1), \ldots, \Sol_{\Omega}(S_k)$
a \emph{Thomas decomposition} of $S$ (with respect to $<$).
\end{defn}

\smallskip

\begin{rem}
Given any (finite) differential system (with coefficients in a
computable differential subfield of $K$), a Thomas decomposition
into simple differential systems can be computed in finitely
many steps by a process which interweaves the algebraic Thomas
algorithm (cf.\ Remark~\ref{rem:algebraicthomasalg}) and
Janet reduction (cf.\ Remark~\ref{rem:janetreduction}).
\end{rem}

\smallskip

\begin{thm}(cf.\ \cite{Robertz2012})\label{thm:simplediff}
Let $S$ as in (\ref{eq:differentialsystem}) be a simple
differential system, $E$ the differential ideal of $R$
generated by ${p_1, \ldots, p_s}$, and let $q$ be the
product of all $\init(p_i)$, $\sep(p_i)$.
Then
\[
E : q^{\infty} := \{ p \in R \mid q^r \cdot p \in E \mbox{ for some } r \in \Z_{\ge 0} \}
\]
is a radical differential ideal, which consists of all
differential polynomials in $R$ vanishing on $\Sol_{\Omega}(S)$.
Given $p \in R$, we have
$p \in E : q^{\infty}$ if and only if $\NF(p, \{ p_1, \ldots, p_s \}) = 0$.
\end{thm}

We refer to \cite{BaechlerGerdtLangeHegermannRobertz2,GerdtThomas,Robertz2012}
for more information and to \cite{ThomasPackages} for an
implementation in Maple.

\section{Elimination}\label{sec:elimination}

We continue to use the same notation as in the previous section.
Our objective is to perform a projection of the solution set of
a differential system onto the space which is addressed by only
certain of the components of the solution tuples. In other words,
we would like to determine all differential consequences of the
given system involving selected differential indeterminates only.

\smallskip

\begin{defn}
Let $B_1$, \ldots, $B_k$ form a partition of the set of
differential indeterminates $\{ u_1, \ldots, u_m \}$.
The \emph{block ranking} on $R$ \emph{with blocks} $B_1$, \ldots, $B_k$
is defined as follows, where $u_{i_1} \in B_{j_1}$, $u_{i_2} \in B_{j_2}$,
$J_1$, $J_2 \in (\Z_{\ge 0})^n$:
\[
\partial^{J_1} u_{i_1} < \partial^{J_2} u_{i_2} \, :\Leftrightarrow \,
j_1 > j_2 \mbox{ or } (j_1 = j_2 \mbox{ and } \partial^{J_1} < \partial^{J_2}).
\]
The comparison of $\partial^{J_1}$ and $\partial^{J_2}$ is
defined by a choice of a degree-reverse lexicographical
ordering as in Example~\ref{ex:degrevlex}. Such a ranking is said to
satisfy $B_1 \gg B_2 \gg \ldots \gg B_k$.
\end{defn}

Using the above notation, for any $1 \le i \le k$, we define
$K\{ B_i, \ldots, B_k \} := K\{ u \mid u \in B_i \cup \ldots \cup B_k \} \subseteq R$.

\smallskip

\begin{prop}(cf.\ \cite{Robertz2012})\label{prop:elimblock}
In the situation of Theorem~\ref{thm:simplediff}, suppose that $<$
is a block ranking with blocks $B_1$, \ldots, $B_k$. For
$1 \le i \le k$, let $E_i$ be the differential ideal of
$K\{ B_i, \ldots, B_k \}$ generated by
$\{ p_1, \ldots, p_s \} \cap K\{ B_i, \ldots, B_k \}$, and let
$q_i$ be the product of all initials and separants of the elements
in this intersection. Then, for every $1 \le i \le k$,
\[
(E : q^{\infty}) \cap K\{ B_i, \ldots, B_k \} = E_i : q_i^{\infty}.
\]
\end{prop}

Hence, computing a Thomas decomposition with respect to a
block ranking enables us to extract generating sets for the
differential consequences satisfied by the projection of the
solution set of the given differential system.

For any differential system $S$ let $S^{=}$ (resp.\ $S^{\neq}$)
denote the set of left hand sides of equations (resp.\
inequations) in $S$.

\smallskip

\begin{cor}(cf.\ \cite{Robertz2012})\label{cor:elimblock}
Let $S$ be a (finite, not necessarily simple) differential system,
and let $S_1$, \ldots, $S_r$ be a Thomas decomposition of $S$
with respect to a block ranking with blocks $B_1$, \ldots, $B_k$.
Let $E$ be the differential ideal of $R$ generated by $S^{=}$
and define the product $q$ of all elements of $S^{\neq}$.
Let $i \in \{ 1, \ldots, k \}$.
For $1 \le j \le r$, let $E_j$ be the differential ideal of
$K\{ B_i, \ldots, B_k \}$ generated by
$S_j^{=} \cap K\{ B_i, \ldots, B_k \}$, and let
$q_j$ be the product of all initials and separants of the elements
in this intersection. Then,
\[
\sqrt{E : q^{\infty}} \cap K\{ B_i, \ldots, B_k \} =
(E_1 : q_1^{\infty}) \cap \ldots \cap (E_r : q_r^{\infty}).
\]
\end{cor}

\section{Systems Theory}\label{sec:systemstheory}

In this section we adapt well-known concepts of nonlinear
control theory to our framework
(cf., e.g.,
\cite{Glad,DiopDiffAlg,FliessGlad,FliessLevineMartinRouchon,Pommaret2001}
and the references therein).

Let $R$ be the differential polynomial ring
$K\{ u_1, \ldots, u_m \}$ in the differential indeterminates
$u_1$, \ldots, $u_m$ with
pairwise commuting derivations $\partial_1$, \ldots, $\partial_n$
as in Subsection~\ref{subsec:diffsimple}.

We assume that a (nonlinear) control system is given by a
simple differential system $S$ as in (\ref{eq:differentialsystem}).
Let $E$ be the differential ideal of $R$
generated by $p_1, \ldots, p_s$, and let $q$ be the
product of all $\init(p_i)$, $\sep(p_i)$.
Let $U := \{ u_1, \ldots, u_m \}$.
(No distinction is made a priori between input, output, state
variables, etc.)

\smallskip

\begin{defn}\label{de:observable}
Let $x \in U$ and
$Y \subseteq U \setminus \{ x \}$.
Then $x$ is \emph{observable with respect to} $Y$
if there exists ${p \in (E : q^{\infty}) \setminus \{ 0 \}}$
such that $p$ is a polynomial in $x$ (not involving any
proper derivative of $x$) with coefficients in
$K\{ Y \}$ and such that neither its leading coefficient
nor $\frac{\partial p}{\partial x}$ is in $E : q^{\infty}$.
\end{defn}

\smallskip

\begin{rem}\label{rem:obervable}
Given the components of a solution to $S$ corresponding to
the variables in $Y$, the implicit function theorem
yields the component corresponding to $x$, if there exists
a polynomial $p$ as in Definition~\ref{de:observable}.
If $S$ is simple with respect to a block ranking satisfying
$U \setminus (Y \cup \{ x \}) \gg \{ x \} \gg Y$ by
Corollary~\ref{cor:elimblock} there exists such a 
$p\in (E:q^\infty)\setminus\{0\}$ if and only if such a 
$p$ exists in $S^=\cap K\{Y,x\}$. Otherwise, one has to
compute a Thomas decomposition of $S$ with respect to such
a ranking and inspect each simple system.
\end{rem}

\smallskip

\begin{defn}\label{de:flat}
Let $Y \subseteq U$. Then
$Y$ is called a \emph{flat output} if
$(E : q^{\infty}) \cap K\{ Y \} = \{ 0 \}$
and every $x \in U \setminus Y$ is observable
with respect to $Y$. The control system given by $S$
is said to be \emph{flat} if a flat output exists.
\end{defn}

\smallskip

\begin{rem}
Whereas, again, computation of a differential Thomas
decomposition of $S$ with respect to a block ranking
satisfying $U \setminus Y \gg Y$ allows to decide
whether $Y$ is a flat output, deciding whether a given
nonlinear control system is flat is a difficult problem
in general.
\end{rem}

Many further applications of the Thomas decomposition method to
systems theory
(e.g., computation of the input-output behaviour from a state space representation,
parameter identification, realization, inversion)
can be realized and will be studied in the future
(cf.\ also \cite{DiopDiffAlg}).

\section{Applications}\label{sec:applications}

In this section we demonstrate the Thomas decomposition method
on two examples.

\smallskip

\begin{exmp}
As an application we consider the model of a continuous stirred-tank
reactor taken from \cite{KwakernaakSivan}.
  
  This model describes a tank with a dissolved material of concentration $c$, which is assumed to be the same everywhere in the tank due to stirring.
  Two input feeds with flow rates $F_1$ and $F_2$ feed this material into the tank with constant concentrations $c_1$ and $c_2$, respectively.
  There exists an outward flow with a flow rate proportional to the square root of the volume $V$ of liquid in the tank.
  The system is modelled by the two differential equations
  \begin{align*}
    \dot{V}(t)                    &= F_1(t) + F_2(t) - k \, \sqrt{V(t)}                                      \\
    \dot{\overline{c(t) \, V(t)}} &= c_1 \, F_1(t) + c_2 \, F_2(t) - c(t) \, k \, \sqrt{V(t)}     \mbox{ ,}
  \end{align*}
  for an experimental constant $k$.
  
  In the following we want to describe the properties of the system in dependence of the constants $c_1$ and $c_2$.
  The Thomas algorithm performs this analysis if we model these contants as functions $c_1(t)$ and $c_2(t)$ satisfying the differential equations $\dot{c}_1(t)=0$ and $\dot{c}_2(t)=0$.
  Additionally, the square root $\sqrt{V(t)}$ of $V(t)$ appears in the equations.
  Our formalism cannot handle roots of functions directly; instead, we introduce $\sqrt{V(t)}$ as new differential indeterminate and substitute $V(t)$ by $\left(\sqrt{V(t)}\right)^2$.
  To exclude trivial cases, we assume $c_1(t)\not=0$, $c_2(t)\not=0$, and $V(t)\not=0$.
  
  We apply our implementation (cf.\ \cite{ThomasPackages}).
  The command \texttt{ComputeRanking} sets the ranking.
  Its first parameter is the list of independent variables and the second parameter is the list of differential indeterminates; a list of lists of differential indeterminates indicates a block ranking.
  We input a derivative of a differential indeterminate as the name of the differential indeterminate indexed by its order; for example \texttt{sV[1]} stands for ${\frac{d}{dt}}\sqrt{V(t)}$.
  The main command \texttt{DifferentialThomasDecomposition} computes a differential Thomas decomposition for a given list of equations and a list of inequations with respect to the fixed ranking.

  In the following we compute a Thomas decomposition of the system using a ranking with $\{F_1,F_2\}\gg \{\sqrt{V},c\}\gg \{c_1,c_2\}$.
  
  \begin{mapleinput}
  \mapleinline{active}{1d}{ivar:=[t]:dvar:=[[F1,F2],[sV,c],[c1,c2]]:
  ComputeRanking(ivar,dvar);
  L:=[ 2*sV[1]*sV[0]-F1[0]-F2[0]+k*sV[0],
    \qquad c[1]*sV[0]\symbol{94}2-c2[0]*F2[0]+c[0]*k*sV[0]
    \qquad\ \ -c1[0]*F1[0]+2*c[0]*sV[1]*sV[0],
    \qquad c1[1], c2[1]]:
  res:=DifferentialThomasDecomposition(
      \quad L,[sV[0],c1[0],c2[0]]);}{}
  \end{mapleinput}
  \mapleresult
  \begin{maplelatex}
  \mapleinline{inert}{2d}{res := [DifferentialSystem, DifferentialSystem, ]}{\[\displaystyle {\it res}\, := \,[{\it DifferentialSystem},{\it DifferentialSystem},\]}
  \mapleinline{inert}{2d}{DifferentialSystem, DifferentialSystem]}{\[\qquad\qquad{\it DifferentialSystem}]\]}
  \end{maplelatex}

  This yields a decomposition consisting of three simple systems.
  We print the first system; for better legibility we have underlined the leaders of equations.  

  \begin{maplegroup}
  \begin{mapleinput}
  \mapleinline{active}{1d}{subs(sV(t)=sqrt(V(t)),
      \quad PrettyPrintDifferentialSystem(res[1]));
  }{}
  \end{mapleinput}
  \mapleresult
  \begin{maplelatex}
  \begin{align}
    [ \left(c2(t) -c1(t)\right) \underline{F1(t)} + \left( {\frac {d}{dt}}c(t)  \right)  \left( \sqrt{V(t)}  \right)^2        &         \nonumber           \\
      +\left(c(t)-c2(t)\right) \left(2\,{\frac {d}{dt}}\sqrt{V(t)} +k\right)\sqrt{V(t)}                                       &= 0,     \label{gen_eq_1}    \\
      \left(c1(t)-c2(t)\right) \underline{F2(t)} +\left( {\frac {d}{dt}}c(t)  \right)  \left( \sqrt{V(t)}  \right)^2          &         \nonumber           \\
      +\left(c(t)-c1(t)\right) \left(2\,{\frac {d}{dt}}\sqrt{V(t)} +k\right)\sqrt{V(t)}                                       &= 0,     \label{gen_eq_2}    \\
      \underline{{\frac {d}{dt}}c1(t)}=0,\ \underline{{\frac {d}{dt}}c2(t)}                                                   &= 0,     \label{gen_param_eq}\\
      \sqrt{V(t)}\not=0,\ c1(t)-c2(t)                                                                                         &\not=0   \nonumber           \\
      c2(t)\not=0,\ c1(t)                                                                                                     &\not=0 ] \nonumber
  \end{align}
  \end{maplelatex}
  \end{maplegroup}
  The equations \eqref{gen_eq_1} and \eqref{gen_eq_2} allow us to solve for $F_1(t)$ and $F_2(t)$ given any $c(t)$ and $V(t)$.
  Thus, we consider $c(t)$ and $V(t)$ a flat output of the system under the additional condition $c_1-c_2\not=0$ on the constants.
  Note that the two equations  $\dot{c}_1(t)=0$ and $\dot{c}_2(t)=0$ in \eqref{gen_param_eq} just model the parameters $c_1$ and $c_2$ as constants.
  
  The other two systems of this decomposition have the additional condition $c_1=c_2$ on the constants.
  This condition prohibits to control the concentration in the tank as both input feeds are equivalent.
  Thus, these systems do not admit $c(t)$ and $V(t)$ as a flat output.
  
  Now we turn our attention to the observability of $\sqrt{V(t)}$ using Remark~\ref{rem:obervable}.
  Therefore, we choose a ranking with $\{\sqrt{V}\}\gg \{c,F_1,F_2\}\gg \{c_1,c_2\}$.
  A Thomas decomposition with this ranking consists of seven systems:
  
  \begin{mapleinput}
  \mapleinline{active}{1d}{ivar:=[t]:dvar:=[[sV],[c,F1,F2],[c1,c2]]:
  ComputeRanking(ivar,dvar);
  res:=DifferentialThomasDecomposition(
    \quad L,[sV[0],c1[0],c2[0]]);}{}
  \end{mapleinput}
  \mapleresult
  \begin{maplelatex}
  \mapleinline{inert}{2d}{}{\[\displaystyle {\it res}\, := \,[{\it DifferentialSystem},{\it DifferentialSystem},\]}
  \mapleinline{inert}{2d}{}{\[\qquad\qquad{\it DifferentialSystem},{\it DifferentialSystem}\]}
  \mapleinline{inert}{2d}{}{\[\qquad\qquad{\it DifferentialSystem},{\it DifferentialSystem},\]}
  \mapleinline{inert}{2d}{}{\[\qquad\qquad{\it DifferentialSystem}]\]}
  \end{maplelatex}
  In the first two systems an equation with $\sqrt{V(t)}$ as leader appears and thus $\sqrt{V(t)}$ is observable.
  For the first system the condition on the parameters for observability is 
  \[
    (c(t)-c_1)F_1(t)+(c(t)-c_2)F_2(t)\not=0\mbox{.}
  \]
  The second system is not physically feasible as it involves negative input feeds due to $F_2(t)\not=0$ and $F_1(t)=-F_2(t)$.
  
  The other five systems include an equation with $\frac {d}{dt}\sqrt{V(t)}$ as leader and thus $\sqrt{V(t)}$ is not observable.
  This follows for the first of these systems as it contains the Wronskian $F_1(t)\dot{F_2}(t)-F_2(t)\dot{F_1}(t)=0$ which makes the inputs linearly dependent.
  In the second system all three concentrations $c(t)$, $c_1$, and $c_2$ are equal and constant.
  In the third system one input feed is zero and the concentration in the tank is equal to the concentration in the other input feed.
  The last two of these systems are not physically feasible because of negative values, as above.
\end{exmp}

\smallskip

\begin{exmp}
  Now we consider a system of partial differential equations from \cite{PommaretQuadrat1997}.
  The system consists of three linear pdes in $\xi_i(\underline{x})=\xi_i(x_1,x_2,x_3)$ for $i=1,2,3$ with a parametric function $a(x_2)$:
  \begin{align}
     0 &= -a(x_2) {\frac{\partial}{\partial x_1}}\xi_1(\underline{x})+{\frac{\partial}{\partial x_1}}\xi_3(\underline{x}) - \left( {\frac{\partial}{\partial {\it x_2}}}a(x_2)  \right) \xi_2(\underline{x}) \nonumber\\
       &+ \frac{1}{2}\,a(x_2) \left(\nabla\cdot\xi(\underline{x})\right)\mbox{,}\label{eq_qp_1}\\
     0 &= -a(x_2) {\frac{\partial}{\partial {\it x_2}}}\xi_1(\underline{x}) +{\frac{\partial}{\partial {\it x_2}}}\xi_3(\underline{x}) \mbox{,}\label{eq_qp_2}\\
     0 &= -a(x_2) {\frac{\partial}{\partial x_3}}\xi_1(\underline{x})+{\frac{\partial}{\partial x_3}}\xi_3(\underline{x}) -\frac{1}{2}\,\left(\nabla\cdot\xi(\underline{x})\right)\mbox{.}\label{eq_qp_3}
  \end{align}
  We model the parameter $a(x_2)$ as function $a(x_1,x_2,x_3)$ in three independent variables satisfying the differential equations 
  \begin{equation}\label{eq_model_a}
    {\frac{\partial}{\partial x_1}}a(x_1,x_2,x_3)=0\mbox{ and }{\frac{\partial}{\partial x_3}}a(x_1,x_2,x_3)=0\mbox{.}
  \end{equation}
  Note that names for derivatives of differential indeterminates now involve multi-indices; for example \texttt{xi1[1,0,0]} stands for ${\frac{\partial}{\partial x_1}}\xi_1(\underline{x})$.
  We apply our implementation with a block ranking $\{\xi_1,\xi_2,\xi_3\}\gg \{a\}$.
  
  \begin{mapleinput}
  \mapleinline{active}{1d}{L:=[ -a[0,0,0]*xi1[1,0,0]+xi3[1,0,0]
    \quad-a[0,1,0]*xi2[0,0,0]+(1/2)*a[0,0,0]*
    \quad(xi1[1,0,0]+xi2[0,1,0]+xi3[0,0,1]),}{}
  \end{mapleinput}
  \begin{mapleinput}
  \mapleinline{active}{1d}{-a[0,0,0]*xi1[0,1,0]+xi3[0,1,0],}{}
  \end{mapleinput}
  \begin{mapleinput}
  \mapleinline{active}{1d}{-a[0,0,0]*xi1[0,0,1]+(1/2)*xi3[0,0,1]
    \quad-(1/2)*xi1[1,0,0]-(1/2)*xi2[0,1,0],
    a[1,0,0],a[0,0,1]]:}{}
  \end{mapleinput}
  \begin{mapleinput}
  \mapleinline{active}{1d}{ivar:=[x1,x2,x3]:dvar:=[[xi1,xi2,xi3],[a]]:
  ComputeRanking(ivar,dvar);}{}
  \end{mapleinput}
  \begin{mapleinput}
  \mapleinline{active}{1d}{res:=DifferentialThomasDecomposition(L,[]);}{}
  \end{mapleinput}
  \mapleresult
  \begin{maplelatex}
  \mapleinline{inert}{2d}{}{\[\displaystyle {\it res}\, := \,[{\it DifferentialSystem},{\it DifferentialSystem},\]}
  \mapleinline{inert}{2d}{}{\[\qquad\qquad{\it DifferentialSystem}]\]}
  \end{maplelatex}
  The Thomas decomposition yields three systems.
  The first system contains no additional condition for $a(x_2)$ except \eqref{eq_model_a}.
  In this generic case \eqref{eq_qp_1}, \eqref{eq_qp_2}, and \eqref{eq_qp_3} are not formally integrable leading to the compatibility condition $\xi_2(x_1,x_2,x_3)=0$.
  The second system includes the additional condition ${\frac{\partial^2}{\partial x_2^2}}a(x_2)=0$, which was calculated in \cite{PommaretQuadrat1997} as the condition to ensure formally integrability of \eqref{eq_qp_1}, \eqref{eq_qp_2}, and \eqref{eq_qp_3}.
  The third system is a special case of the second system with the new condition $a(x_2)=0$.
\end{exmp}


\begin{ack}                               
The authors would like to thank A.~Quadrat for discussion and
valuable comments related to this paper.
The first author was partially supported by
Schwerpunkt SPP 1489 of the Deutsche Forschungsgemeinschaft.
\end{ack}


\end{document}